\documentclass[11pt]{amsart}
\usepackage[paperheight=11in,paperwidth=8.5in,left=1in,top=1.25in,right=1in,bottom=1.25in,]{geometry}
\usepackage{amsmath,amsthm,comment}

\usepackage{enumitem}
\usepackage{relsize}
\usepackage{rotating}
\usepackage[table]{xcolor}

\usepackage[super]{nth}
\usepackage[none]{hyphenat}
\usepackage{hyperref}
\hypersetup{
	colorlinks,
	citecolor=black,
	filecolor=black,
	linkcolor=black,
	urlcolor=black
}

\usepackage[capitalize, nameinlink, noabbrev, nosort]{cleveref}

\usepackage{graphicx, amssymb}
\usepackage{epstopdf}
\usepackage{youngtab}
\usepackage[boxsize=.35 em]{ytableau}
\usepackage{comment}
\usepackage{setspace}
\usepackage{multicol}

\usepackage{tikz}
\usetikzlibrary{arrows,calc,intersections,matrix,positioning,through}
\usepackage{tikz-cd}
\usepackage{comment}
\usepackage{diagbox}
\usepackage{blkarray}
\tikzset{commutative diagrams/.cd,every label/.append style = {font = \normalsize}}

\newcommand{\lr}[1]{\langle #1 \rangle}
\newcommand{\llrr}[1]{\langle\!\langle #1 \rangle\!\rangle}

\newtheorem{theorem}{Theorem}

\newtheorem{proposition}[theorem]{Proposition}

\theoremstyle{definition}
\newtheorem{definition}[theorem]{Definition}
\newtheorem{example}[theorem]{Example}

\newtheorem{notation}[theorem]{Notation}

\DeclareMathOperator{\Gr}{Gr}
\DeclareMathOperator{\Mat}{Mat}

\DeclareMathOperator{\refl}{refl}
\DeclareMathOperator{\cyc}{cyc}

\DeclareMathOperator{\pre}{pre}

\DeclareMathOperator{\Irr}{\xx}

\newcommand{\Grk}{\Gr_{k,n}^{\scriptscriptstyle \ge 0}}

\def\AA{\mathcal{A}_{n,k,m}(Z)}
\newcommand{\Ank}{\mathcal{A}_{n, k, 4}({Z})}

\def\bcfw{\bowtie}

\def\4biddenprop{4-coindependent}

\newcommand{\gt}[1]{Z_{#1}}
\newcommand{\gto}[1]{Z_{#1}^\circ}

\newcommand{\rzeta}{\bar{\zeta}}
\newcommand{\ralpha}{\bar{\alpha}}
\newcommand{\rbeta}{\bar{\beta}}
\newcommand{\rgamma}{\bar{\gamma}}
\newcommand{\rdelta}{\bar{\delta}}
\newcommand{\repsilon}{\bar{\varepsilon}}

\newcommand{\czeta}{{\zeta}}

\newcommand{\rPsi}{\overline{\Psi}}

\newcommand{\C}{\mathbb{C}}
\newcommand{\CC}{\mathbb{C}}
\newcommand{\A}{\mathcal{A}}

\def\Xcal{\mathcal{X}}
\def\Acal{\mathcal{A}}

\def\Fcal{\mathcal{F}}

\def\R{\mathbb{R}}

\newcommand{\xx}{\mathbf{x}}

\usepackage{lipsum}
\usepackage{wrapfig}
\usepackage[rightcaption]{sidecap}

\begin{document}
\begin{abstract}
In 2005, Britto, Cachazo, Feng and Witten gave a recurrence (now known as the BCFW recurrence)
for computing scattering amplitudes in N=4 super Yang Mills theory. Arkani-Hamed and Trnka
subsequently introduced the amplituhedron to give a geometric interpretation of the BCFW recurrence. Arkani-Hamed and Trnka conjectured that each way of iterating the BCFW recurrence gives a ``triangulation'' or ``tiling'' of the m=4 amplituhedron. In this article we prove the BCFW tiling conjecture of Arkani-Hamed and Trnka. We also prove the cluster adjacency conjecture for BCFW tiles of the amplituhedron, which says that
facets of tiles are cut out by collections of compatible cluster variables for the Grassmannian $\Gr_{4,n}$.
Moreover we show that each BCFW tile is the subset of the Grassmannian where certain cluster variables have particular~signs.

\medskip
\noindent
This article was published on PNAS
(\href{https://doi.org/10.1073/pnas.2408572122}{doi.org/10.1073/pnas.2408572122}) as a research announcement of a full-length paper by the same authors (\href{https://doi.org/10.48550/arXiv.2310.17727}{doi.org/10.48550/arXiv.2310.17727}).
\end{abstract}

\title[BCFW tilings and cluster adjacency for the amplituhedron]{BCFW tilings and cluster adjacency \\ for the amplituhedron}

\date{}

\author[C. Even-Zohar]{Chaim Even-Zohar}
\address{Faculty of Mathematics, Technion, Haifa, Israel}
\email{chaime@technion.ac.il}

\author[T. Lakrec]{Tsviqa Lakrec}
\address{Section de math\'ematiques, Universit\'e de Gen\`eve, Gen\`eve, Switzerland}
\email{tsviqa@gmail.com}

\author[M. Parisi]{Matteo Parisi}
\address{Institute for Advanced Study, Princeton, NJ; CMSA, Harvard University, Cambridge, MA}
\email{mparisi@cmsa.fas.harvard.edu}
	
\author[M. Sherman-Bennett]{Melissa Sherman-Bennett}
\address{Department of Mathematics, MIT, Cambridge, MA}
\email{msherben@mit.edu}
	
\author[R. Tessler]{Ran Tessler}
\address{Department of Mathematics, Weizmann Institute of Science, Israel}
\email{ran.tessler@weizmann.ac.il}
	
\author[L. Williams]{Lauren Williams}
\address{Department of Mathematics, Harvard University, Cambridge, MA}
\email{williams@math.harvard.edu}

\maketitle


\begin{figure}
\centering  \includegraphics[width=0.47\linewidth]{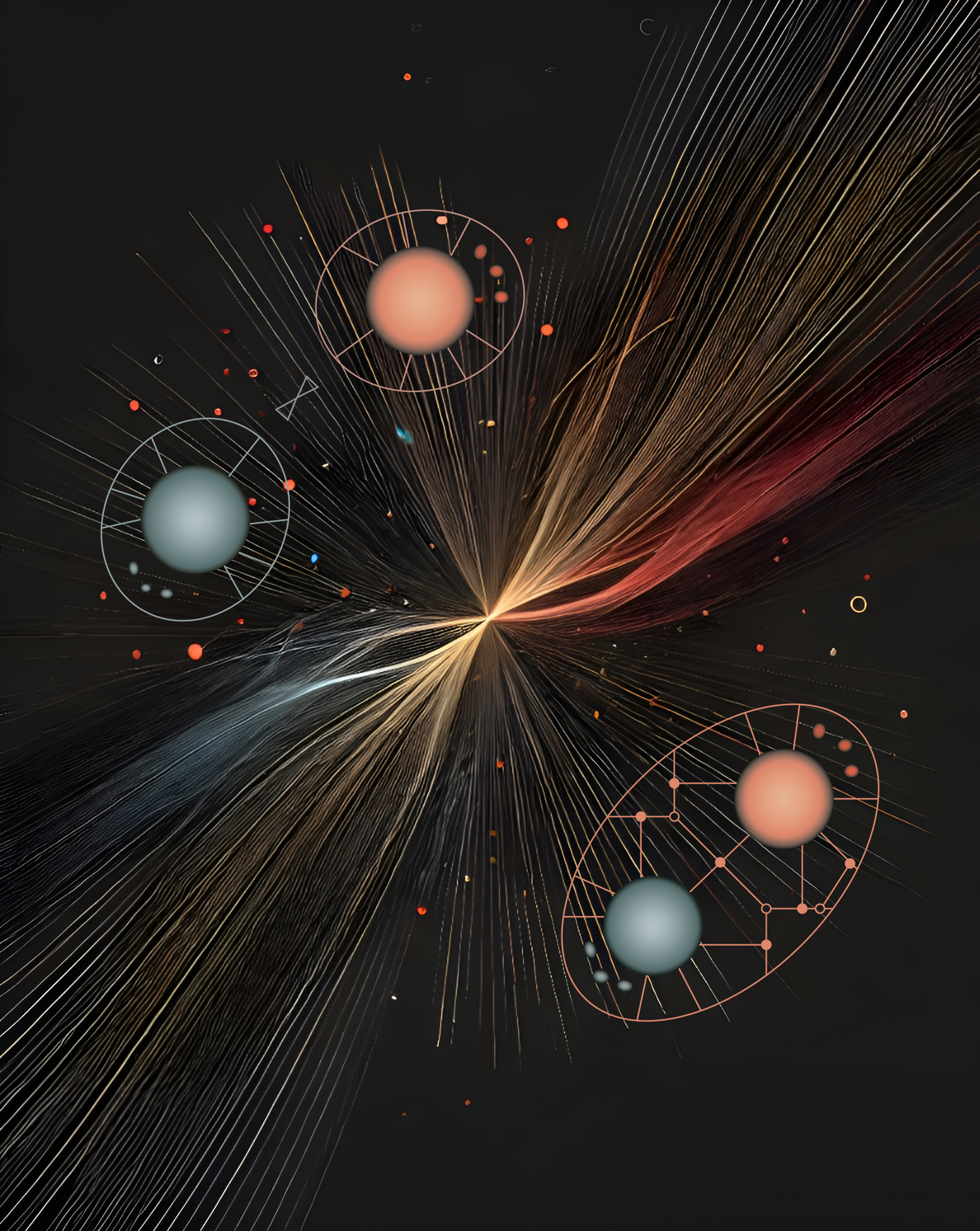}
\caption{Artistic illustration by Annabel Ma}
\label{fig:enter-label}
\end{figure}

{T}he (tree) \emph{amplituhedron} $\AA$ is the image of the
positive Grassmannian $\Gr_{k,n}^{\scriptscriptstyle\geq 0}$ under the
\emph{amplituhedron map} $\tilde{Z}: \Gr_{k,n}^{\scriptscriptstyle\geq 0} \to \Gr_{k,k+m}$.
It was introduced by
Arkani-Hamed and Trnka \cite{arkani-hamed_trnka} in order  to give a
geometric  interpretation of
\emph{scattering amplitudes} in $\mathcal{N}=4$ super Yang Mills theory (SYM):
in particular, they conjectured one can
 compute $\mathcal{N}=4$ SYM scattering amplitudes  by
`tiling' the $m=4$ amplituhedron $\Ank$
--- that is, decomposing the amplituhedron into
`tiles' --- and summing the `volumes' of the tiles.
While the case $m=4$ is most important for physics, the amplituhedron
is defined for any positive $n,k,m$ with $k+m \leq n$, and
has a very rich geometric and combinatorial
structure.  It generalizes cyclic polytopes (when $k=1$),
cyclic hyperplane arrangements \cite{karpwilliams}
(when $m=1$), and the positive Grassmannian (when $k=n-m$), and it
is connected to
the hypersimplex and the positive tropical Grassmanian \cite{LPW,PSW} (when $m=2$).
The present work proves the \emph{BCFW tiling conjecture} and the \emph{cluster adjacency conjecture} for BCFW tiles, addressing two longstanding problems for the amplituhedron.

The \emph{cluster adjacency conjecture} says that
facets of tiles are cut out by collections of compatible cluster variables.
This was motivated by the observation that cluster algebras can be used to  describe singularities of scattering
amplitudes in $\mathcal{N}=4$ SYM \cite{Golden:2013xva}. In particular, \cite{Drummond:2017ssj, Drummond:2018dfd} conjectured that
the terms in tree-level amplitudes
coming from the BCFW recursions are rational functions
whose poles correspond to compatible cluster variables of the Grassmannian $\Gr_{4,n}$, see also \cite{Mago:2019waa}.  The cluster adjacency
conjecture, formulated for the $m=2$ and
$m=4$ amplituhedron in
\cite{Lukowski:2019sxw} and
 \cite{Gurdougan2020ClusterPI}, was subsequently proved for all tiles of
the $m=2$ amplituhedron  \cite{PSW}, but remained open in the $m=4$ case,
the case of most relevance to physics.

The \emph{BCFW tiling conjecture} says that any way of iterating the
BCFW recurrence \cite{BCFW} gives rise to a collection of cells (``BCFW cells'') in the positive Grassmannian whose images tile the
$m=4$ amplituhedron $\Ank$ (see \cref{def:tiling} for a precise definition).  This conjecture arose side-by-side with the definition of the amplituhedron -- indeed, the goal of Arkani-Hamed and Trnka, realized in \cite{arkani-hamed_trnka}, was to find a geometric object which could be decomposed into pieces coming from BCFW cells.
BCFW-like tilings of the $m=1$ and $m=2$ amplituhedron
were proved in \cite{karpwilliams}
and \cite{BaoHe}, building on \cite{AHTT}
and \cite{karp2020decompositions}.
A step towards the BCFW tiling conjecture was made in \cite{even2021amplituhedron}, where the authors built on work
of \cite{karp2020decompositions} to show that the `standard' way of performing the BCFW recursion gives a tiling for the $m=4$ amplituhedron.

\subsection*{Main results.}
In this paper we extend and generalize the results of \cite{PSW} and
 \cite{even2021amplituhedron}
 to give a very
complete picture of the $m=4$
amplituhedron. We show that arbitrary \emph{BCFW cells} give tiles (\cref{thm:bcfw_tiles}) and that they satisfy the cluster adjacency conjecture (\cref{thm:clusteradjacency}).
We strengthen the connection with cluster algebras by associating to each BCFW tile a collection of compatible cluster variables for $\Gr_{4,n}$ (\cref{def:coord_cluster_vars}), which we use to describe the tile as a semialgebraic set in $\Gr_{k,k+4}$ (\cref{th:sign_description}).
For `standard' BCFW tiles, one can also give a non-recursive description of these cluster variables \cite[Theorem 8.4]{Even-Zohar:2023del} and the underlying quiver (\cref{thm:quiver}).
Finally, we use these results to prove the BCFW tiling conjecture for the $m=4$ amplituhedron (\cref{thm:BCFWtiling}).

\subsection*{Further motivation.}
From the point of view of cluster algebras,
the study of tiles for the amplituhedron $\mathcal{A}_{n,k,m}$ is useful
because it is closely related to the
cluster structure on the Grass- mannian $\Gr_{m,n}$, as was shown for
$m=2$ in
\cite{PSW} and as this paper demonstrates for $m=4$.  In particular, for $m=4$,
the \emph{BCFW product} (\cref{def:butterfly}) used to recursively build tiles has a cluster quasi-homomorphism
counterpart called \emph{product promotion} (\Cref{def:product_promotion}),
that can be used to recursively construct
cluster variables and seeds in $\Gr_{4,n}$ (\cref{thm:promotion2}). Moreover, this would also have potential implications for loop amplitudes in $\mathcal{N}=4$ SYM in connection with the \emph{cluster bootstrap} program (see \cite{Caron-Huot:2020bkp} for a review).

In the closely related field of \emph{total positivity}, one prototypical
problem is to give an efficient characterization of
the `positive part' of a space as the subset
where a collection of functions take on positive values \cite{fz-intel}. A minimal collection of functions whose positivity cuts out the positive part is called a `positivity test'.
For example,  for any cluster $\xx$ for $\Gr_{k,n}$ \cite{scott}, the \emph{positive Grassmannian} $\Gr_{k,n}^{>0}$ can be described as the region in  $\Gr_{k,n}$ where all the cluster variables of $\xx$ are positive,
so each $\xx$ provides a positivity test. 

We think of \cref{th:sign_description}
as a `positivity test'
for membership in a BCFW tile of the amplituhedron.
See \cite[Theorem 6.8]{PSW} for an analogous result for $m=2$,
and \cite[Conjecture 7.17]{Even-Zohar:2023del} for some conjectures for general $m$.

From the point of view of discrete geometry,
one can think of about tiles as a generalization of polytopes in the Grassmannian.
In particular, the positivity tests for the positive Grassmannian and BCFW tiles
can be thought of as
analogues of the hyperplane description of polytopes.
Finally, it is expected that
   tiles are \emph{positive geometries} \cite{Arkani-Hamed:2017tmz}.


\section*{Background}

Let $[n]$ denote $\{1,\dots,n\}$, and $\binom{[n]}{k}$ denote the set of all $k$-element
subsets of $[n]$.

 \subsection*{The (positive) Grassmannian}\label{sec:posGrass}

The \emph{Grassmannian} $\Gr_{k,n}$
is the space of all $k$-dimensional subspaces of $\mathbb{C}^n$. We denote its real points by $\Gr_{k,n}(\mathbb{R})$.
We can
 represent a point $V \in
\Gr_{k,n}$ as the row-span
of
a full-rank $k\times n$ matrix $C$.
For $I=\{i_1 < \dots < i_k\} \in \binom{[n]}{k}$, we let $\lr{I}_V=\lr{i_1\,i_2\,\dots\,i_k}_V$ be the $k\times k$ minor of $C$ using the columns $I$.
The $\lr{I}_V$
are called the {\itshape Pl\"{u}cker coordinates} of $V$, and are independent of the choice of matrix
representative $C$ (up to common rescaling). The \emph{Pl\"ucker embedding}
$V \mapsto \{\lr{I}_V\}_{I\in \binom{[n]}{k}}$
embeds $\Gr_{k,n}$ into
projective space\footnote{
We will sometimes abuse notation and identify $C$ with its row-span;
we will also drop the subscript $V$ on Pl\"ucker coordinates when it does not cause confusion.}.
If $C$ has columns $v_1, \dots, v_n$, we may also identify $\lr{i_1\,i_2\,\dots\,i_k}$ with $v_{i_1} \wedge v_{i_2} \wedge \dots \wedge v_{i_k}$, hence e.g. $\lr{i_1\,i_2\,\dots\,i_k}=- \lr{i_2\,i_1\,\dots\,i_k}$.
We occasionally will also work with $\Gr_{k,N}$, the Grassmannian of $k$-planes in a vector space with basis indexed by $N\subset [n]$.

\begin{definition}[Positive Grassmannian \cite{lusztig, postnikov}]\label{def:positroid}
We say that $V\in \Gr_{k,n}(\R)$ is \emph{totally nonnegative}
     if (up to a global change of sign)
       $\lr{I}_V \geq 0$ for all $I \in \binom{[n]}{k}$.
Similarly, $V$ is \emph{totally positive} if $\lr{I}_V >0$ for all $I
      \in \binom{[n]}{k}$.
We let $\Grk$ and $\Gr_{k,n}^{>0}$ denote the set of
totally nonnegative and totally positive elements of $\Gr_{k,n}$, respectively.
$\Grk$ is called the \emph{totally nonnegative}  \emph{Grassmannian}, or
       sometimes just the \emph{positive Grassmannian}.
\end{definition}


If we partition $\Grk$ into strata based on which Pl\"ucker coordinates are strictly
positive and which are $0$, we obtain a cell decomposition of $\Grk$
into \emph{positroid cells} \cite{postnikov}.
Each positroid cell $S$ gives rise to a matroid $\mathcal{M}$, whose bases are precisely
the $k$-element subsets $I$ such that the Pl\"ucker coordinate
$\lr{I}$ does not vanish on $S$; $\mathcal{M}$
is called a \emph{positroid}.
There are many ways to index positroid cells in $\Grk$ \cite{postnikov}, such as \emph{plabic graphs}.

\begin{definition}
	Let $G$ be a \emph{plabic graph}\footnote{We will always assume that plabic graphs are \emph{reduced} \cite[Definition 12.5]{postnikov}.}, i.e. a planar graph embedded in a disk, with boundary vertices $1, 2, \dots, n$ on the boundary of the disk and with internal vertices colored black or white. A \emph{perfect orientation} $\mathcal{O}$ of $G$ is an orientation of the edges so that each black vertex has a unique outgoing edge and each white vertex has a unique incoming edge. Each perfect orientation has the same number of sources, which are boundary vertices. The positroid associated to $G$ is the collection $\mathcal{M}:=\{I: I \text{ the source set of a perfect orientation of }G \}$.
\end{definition}

Both $\Gr_{k,n}$ and $\Grk$ admit the following set of operations,
which will be useful to us.

\begin{definition}[Operations on the Grassmannian]\label{def:opGr} We define the following maps on $\Mat_{k,n}$, which descend to maps (denoted in the same way) on $\Gr_{k,n}$ and $\Grk$:
\begin{itemize}[leftmargin=20pt]
\item (cyclic shift) The map $\cyc: \Mat_{k, n} \to \Mat_{k,n}$ sends  $v_1 \mapsto (-1)^{k-1}v_{n}$ and $v_i \mapsto v_{i-1}, 2 \leq i \leq n$. 
 \item  (reflection) The map $\refl:\Mat_{k, n} \to \Mat_{k,n}$ sends $v_i \mapsto v_{n+1-i}$ and rescales the top row by~$(-1)^{\binom{k}{2}}$.
 \item (zero column) The map $\pre_i:\Mat_{k, [n] \setminus \{i\}} \to \Mat_{k,n}$ adds a zero column at $i$. 
\end{itemize}
\end{definition}

\subsection*{The amplituhedron}

Building on \cite{abcgpt,hodges},
Arkani-Hamed and Trnka
\cite{arkani-hamed_trnka}
intro- duced
the \emph{(tree) amplituhedron}, which is
the image of the positive Grassmannian under a positive linear map.
Let $\Mat_{n,p}^{>0}$ denote the set of $n\times p$ matrices whose maximal minors
are positive.

\begin{definition}[Amplituhedron]\label{defn_amplituhedron}
Let $Z\in \Mat_{n,k+m}^{>0}$, where $k+m \leq n$. 
    The \emph{amplituhedron map}
$\tilde{Z}:\Gr_{k,n}^{\ge 0} \to \Gr_{k,k+m}$
        is defined by
        $\tilde{Z}(C):=CZ$,
    where as usual we identify the matrices $C, CZ$ with their rowspans.
        The \emph{amplituhedron} $\AA \subset \Gr_{k,k+m}$ is the image
$\tilde{Z}(\Gr_{k,n}^{\ge 0})$.
\end{definition}

The amplituhedron is $km$-dimensional, so it is full-dimensional in $\Gr_{k,k+m}$. We will be concerned with decompositions of $\AA$ into full-dimensional images of particular positroid cells.

\begin{definition}[Tiles]
\label{defn_tile}
	Fix $k, n, m$ with $k+m \leq n$ and choose
	$Z\in \Mat_{n,k+m}^{>0}$.
	Given a positroid cell $S$ of
	$\Gr_{k,n}^{\ge 0}$, we let
	$\gto{S} := \tilde{Z}(S)$ and
	$\gt{S}: = \overline{
		\tilde{Z}(S)} = \tilde{Z}(\overline{S})$.
	We call
	$\gt{S}$  and $\gto{S}$
	a \emph{tile}  and an \emph{open tile}
	for $\AA$
	if $\dim(S) =km$ and
	$\tilde{Z}$ is injective on $S$.
\end{definition}

\begin{definition}[Tilings]\label{def:tiling}
	A \emph{tiling}
                of $\AA$ is
                a collection
		$\{Z_{S} \ \vert \ S\in \mathcal{C}\}$
                 of tiles, such that their union equals $\AA$ and all open tiles in $\{\gto{S} \ \vert \ S\in \mathcal{C}\}$ are pairwise disjoint.
\end{definition}

There is a natural notion of  \emph{facet} of a tile,
generalizing the notion of facet of a polytope.

\begin{definition}[Facet of a cell and a tile]\label{def:facet2}
Given two positroid cells $S'$ and $S$, we say that
$S'$ is a \emph{facet} of $S$ if
$S' \subset \partial{S}$ and $S'$ has codimension $1$ in $\overline{S}$.
If $S'$ is a facet of $S$ and $Z_S$ is a tile of $\AA$, we say that $Z_{S'}$ is a
	\emph{facet}
 of $Z_{S}$ if
         $\gt{S'} \subset \partial \gt{S}$ and has codimension 1 in $\gt{S}$.
\end{definition}

In order to describe tiles and their facets, we will use functions that are adapted to the amplitu- hedron map.

\begin{definition}[Twistor coordinates]\label{def:tw_coords}
Fix
  $Z \in \Mat_{n,k+m}^{>0}$
   with rows  $Z_1,\dots, Z_n \in \R^{k+m}$.
   Given  $Y \in \Gr_{k,k+m}$
   with rows $y_1,\dots,y_k$,
   and $\{i_1,\dots, i_m\} \subset [n]$,
   we define the \emph{twistor coordinate} $
	\llrr{{i_1} {i_2} \cdots {i_m}}$
        to be  the determinant of the
        matrix with rows
        $y_1, \dots,y_k, Z_{i_1}, \dots, Z_{i_m}$.
\end{definition}

Note that the twistor coordinates are defined only up to a common scalar multiple. An element of $\Gr_{k, k+m}$ is uniquely determined by its twistor coordinates \cite{karpwilliams}. Moreover, $\Gr_{k,k+m}$ can be embedded into $\Gr_{m,n}$ so that the twistor coordinate $\llrr{i_1 \dots i_m}$ is the pullback of the Pl\"ucker coordinate $\lr{i_1, \dots, i_m}$ in $\Gr_{m,n}$.
\begin{definition}\label{def:functionary}
We refer to a homogeneous polynomial in twistor coordinates as a \emph{functionary}. For $S \subseteq \Grk$, we say a functionary $F$ has a definite sign $s \in \{\pm 1\}$ on $\gto{S}$ if for all $Z\in \Mat_{n,k+4}^{>0}$ and for all $Y \in \gto{S}$, $F(Y)$ has sign $s$.
\end{definition}

We will use functionaries to describes BCFW tiles.

\subsection*{Cluster Algebras}

Cluster algebras were introduced by Fomin and Zelevinsky in \cite{FZ1}, motivated by the study of total positivity; see \cite{FWZ} for
an introduction.
We quickly recall notation for cluster algebras of geometric type from quivers.

A \emph{quiver} $Q$ is a finite directed graph.
We require that quivers do not have any oriented cycles of length $1$ or $2$.

\begin{definition}
\label{def:seed0}
Choose $s\geq r$ positive integers.
Let $\Fcal$ be a field
of rational functions
in $r$ independent
variables
over
$\CC(x_{r+1},\dots,x_s)$.
A \emph{labeled seed} in~$\Fcal$ is
a pair $(\xx, Q)$, where $\xx = (x_1, \dots, x_s)$ forms a free generating
set for $\Fcal$ and $Q$ is a quiver on $[s]$ where vertices $[r]$ are called \emph{mutable} and vertices $r+1, \dots, s$ are called \emph{frozen}.
The tuple $\xx$ is a \emph{cluster}. The functions $\{x_1, \dots, x_r\}$ are \emph{mutable cluster variables} and the functions $c=\{x_{r+1},\dots,x_s\}$ are \emph{frozen cluster variables}.
\end{definition}

An operation called \emph{mutation} produces new seeds \cite{FZ1}.

\begin{definition}
[{Quiver Mutation}]
Let $Q$ be a quiver and let $k$ be a vertex.
The
\emph{mutated quiver} $\mu_k(Q)$ has the same vertex set as $Q$
and its set of arrows is obtained as follows:
\begin{enumerate}[leftmargin=20pt]
\item for each subquiver $i \to k \to j$, add a new arrow $i \to j$;
\item reverse all arrows with source or target $k$;
\item remove all $2$-cycles.
\end{enumerate}
\end{definition}

\begin{definition}[Seed mutation]
    Let $(\xx, Q)$ be a seed in $\mathcal{F}$ and let $k$ be a mutable vertex of $Q$. The mutated seed $\mu_{k}(\xx, Q)=(\xx', Q')$ is another seed in $\Fcal$, with $Q' = \mu_k(Q)$ and $\xx'= \xx\setminus \{x_k\} \cup \{x_k'\}$, where the new cluster variable $x_k'$ is determined by the \emph{exchange relation}
\begin{equation*}
\label{exchange relation0}
x'_k\ x_k =
 \ \prod_{i \rightarrow k \text{ in }Q} x_{i}
+ \ \prod_{i \leftarrow k \text{ in }Q} x_{i} .
\end{equation*}
\end{definition}



\begin{definition}
[{Cluster algebra}]
\label{def:cluster-algebra0}
Given a seed $(\xx, Q)$ in $\mathcal{F}$, we denote as $\Xcal$
the union of all mutable variables of all the seeds obtained from $(\xx, Q)$ by any sequence of mutations.
Let $\CC[c^{\pm 1}]$ be the \emph{ground ring} consisting
of Laurent polynomials in the frozen variables. The
\emph{cluster algebra} $\Acal(\xx, Q)$ is the $\CC[c^{\pm 1}]$-subalgebra of $\Fcal$
generated by all mutable variables,
with coefficients which are Laurent polynomials
in the frozen variables: $\Acal(\xx, Q) = \CC[c^{\pm 1}] [\Xcal]$.
A subset $C$ of cluster variables in $\Xcal$ are \emph{compatible} if there exists a seed $(\xx', Q')$ obtained from $(\xx, Q)$ by a sequence of mutations such that $C$ is a subset of the cluster $\xx'$.
\end{definition}

\begin{figure}
\centering
\includegraphics[width=.5\linewidth]{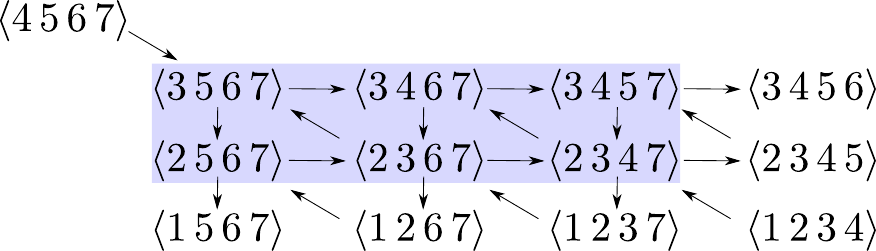}
\caption{The rectangle seed $\Sigma_{4,7}$. Mutable variables are in the colored box.}
\label{fig:gr_quiver}
\end{figure}

The coordinate ring of the Grassmannian $\Gr_{k,n}$
is a cluster algebra \cite{scott}. One may take the \emph{rectangles seed} $\Sigma_{k,n}$ as the initial seed \cite{FW6}, see \cref{fig:gr_quiver}.

\begin{theorem}[\cite{scott}]\label{thm:scott}
Let ${\Gr}_{k,n}^\circ$ be the open subset of the Grassmannian where the frozen variables
don't vanish. Then the coordinate ring $\C[\widehat{\Gr}_{k,n}^\circ]$ of the affine cone over ${\Gr}_{k,n}^\circ$ is the cluster algebra
 $\A(\Sigma_{k,n})$.
\end{theorem}

Moreover, the operations on the Grassmannian $\cyc, \refl, \pre$ in \cref{def:opGr} induce maps on $\CC[\widehat{\Gr}^\circ_{k,N}]$ which are compatible
with the cluster structure.


\begin{proposition}\label{prop:pre_cluster}
       The pullbacks $\cyc^*, \refl^*: \CC[\widehat{\Gr}^\circ_{k,n}] \to \CC[\widehat{\Gr}^\circ_{k,n}], \pre_i^*: \C[\widehat{\Gr}^\circ_{k, n}] \to \C[\widehat{\Gr}^\circ_{k,[n] \setminus \{i\}}]$ take cluster variables to cluster variables and preserve compatibility and exchange relations.
\end{proposition}

\section*{Results}
Our first main result is that a class of cells called \emph{BCFW cells} give rise to tiles for $\Ank$. We will build BCFW cells recursively using the \emph{BCFW product}. Let us first introduce some notation we will use throughout this section.

\begin{notation}\label{not:LR_cluster}
Choose integers $1\leq a<b<c<d<n$ with $a,b$ and $c,d,n$ consecutive.
	Let\footnote{
Note that we will overload the notation and let $n$
	index an element of a vector space basis for different
	vector spaces; however, in what follows, the meaning should
	be clear from context.} $N_L = \{1,2,\dots,a,b,n\}, N_R = \{b, \dots, c, d, n\}$ and $B=(a,b,c,d,n)$. Also fix $k \leq n$ and two nonnegative integers $k_L \leq |N_L|$ and $ k_R\leq |N_R|$ such that $k_L + k_R +1=k$.
 Note that, for any set of indices $N \subset [n]$, our results hold with $N$ instead of $[n]$, by replacing $1$ and $n$ in the definition with the smallest and largest elements of $N$, respectively.
\end{notation}

\begin{definition}[BCFW product]\label{def:butterfly}
Using \cref{not:LR_cluster}, let $S_L \subseteq \Gr^{\geq 0}_{k_L,N_L},\ S_R \subseteq \Gr^{\geq 0}_{k_R,N_R}$ be positroid cells and $G_L, G_R$ be the respective plabic graphs. The \emph{BCFW product} of $S_L$ and $S_R$ is the positroid cell $S_L \bcfw S_R \subseteq \Gr_{k,n}$ corresponding to the plabic graph in the right-hand side of \Cref{fig:butterfly}.

\begin{figure}[h]
\centering
\includegraphics[width=\linewidth]{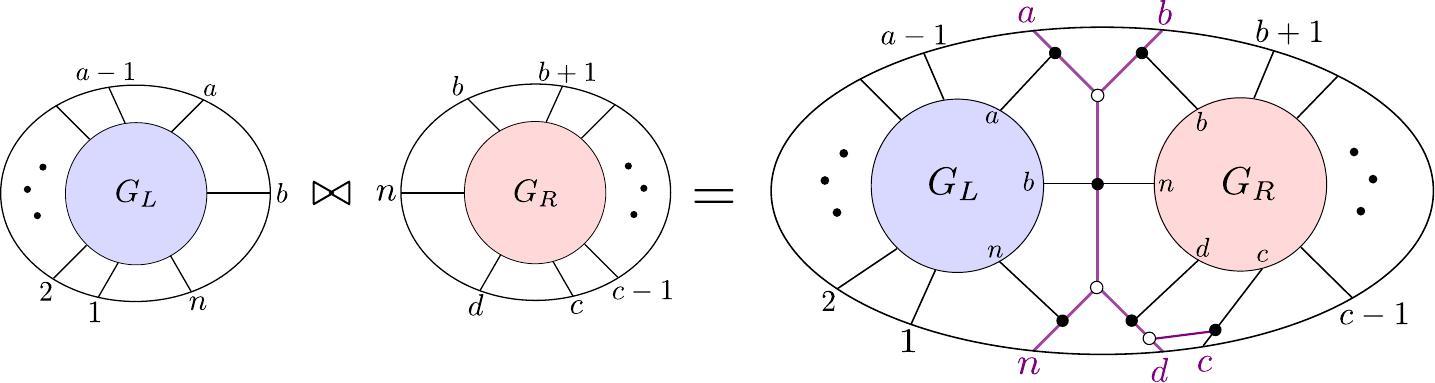}
\caption{The BCFW product $S_L \bcfw S_R$ of $S_L$ and $S_R$ in terms of their plabic graphs.}
\label{fig:butterfly}
\end{figure}

\end{definition}

We now define the family of \emph{BCFW cells} to be the set of positroid cells which is closed under the operations in \Cref{def:opGr,def:butterfly}.

\begin{definition}[BCFW cells]\label{def:BCFW_cell}
The set of \emph{BCFW cells} is defined recursively. For $k=0$, the trivial cell $\Gr^{\scriptscriptstyle>0}_{0,n}$ is a BCFW cell. If $S$ is a BCFW cell, so are $\cyc S$, $\refl S$ and $\pre_i S$. If $S_L,S_R$ are BCFW cells, so is their BCFW product $S_L \bcfw S_R$.
\end{definition}

\begin{figure}
\centering
\vspace{0.5cm}
\includegraphics[width=0.4\linewidth]{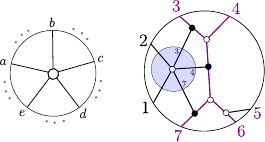}
 \caption{ \label{fig:bcfwcellsk1} Plabic graphs of
   a BCFW cell in $\Gr^{\geq 0}_{1,n}$ (left) and in
   $\Gr^{\geq 0}_{2,7}$ (right).} \label{fig:bcfw_tile}
\end{figure}

\begin{example}\label{ex:bcw_tile}
For $k=1$, the BCFW cells in $\Gr_{1,n}^{\geq 0}$ have plabic graphs as in \cref{fig:bcfw_tile} (left). The associated positroid has bases $\{a\},\{b\},\{c\},\{d\},\{e\}$.
    In \cref{fig:bcfw_tile} (right), $S_{ex} \subset \Gr^{\geq 0}_{2,7}$ is obtained as $S_L \bowtie S_R$, with $S_L, S_R$ BCFW cells in $\Gr^{\geq 0}_{1,N_L}, \Gr^{\geq 0}_{0,N_R}$ respectively, with  $N_L=\{1,2,3,4,7\}, N_R=\{4,5,6,7\}$ and $B=(3,4,5,6,7)$.
\end{example}


Our first main result shows that images of BCFW cells are in fact tiles in $\Ank$.

\begin{theorem}[BCFW tiles]\label{thm:bcfw_tiles}
The amplituhedron map is injective on each BCFW cell.
That is,
	the closure
	${\gt{S} :=
	\overline{\tilde{Z}(S)}}$ of
	the image of a BCFW cell $S$
 is a \emph{tile}, which we
	refer to as a \emph{BCFW tile}.
\end{theorem}

A key ingredient to prove \cref{thm:bcfw_tiles} is inverting the amplituhedron map on BCFW tiles \cite[Theorem 7.7]{Even-Zohar:2023del} by using \emph{product promotion} -- an operation which interacts nicely both with the cluster structure on the Grassmannian and with the BCFW product.

\begin{definition}\label{def:product_promotion} Given 5 vectors $v_i,v_j,v_r,v_s,v_q$, we set $(ij)\cap (rsq):=v_i \lr{j \, r \, s\, q}-v_j \lr{i \, r \, s\, q}= - v_r \lr{i \, j \, s\, q}+v_s \lr{i \, j \, r\, q}-v_q \lr{i \, j \, r\, s} $, which is in the intersection of the subspace spanned by $v_i, v_j$ and that spanned by $v_r, v_s, v_q$.
Using \cref{not:LR_cluster}, \emph{product promotion} is the homomorphism
	$$\Psi_{B} = \Psi: \C(\widehat{\Gr}_{4,N_L})\times \C(\widehat{\Gr}_{4,N_R}) \to \C(\widehat{\Gr}_{4,n}),$$
induced
by the following substitution:
\begin{equation*}
\text{on $\widehat{\Gr}_{4,N_L}$: } v_b \;\mapsto\; \frac{(ba)\cap (cdn)}{\lr{a\,c\,d\,n}},
\end{equation*}
\begin{equation*}
\text{on $\widehat{\Gr}_{4,N_R}$: } v_n \;\mapsto\;
\frac{(ba)\cap (cdn)}{\lr{a\,b\,c\,d}}
\label{eq:promotionvectors2}, \, v_d \;\mapsto\; \frac{(dc)\cap (abn)}{\lr{a\,b\,c\,n}}.
\end{equation*}
\end{definition}

 We show that $\Psi$ is a \emph{quasi-homomorphism} from the cluster algebra\footnote{$\C[\widehat{\Gr}_{4,N_L}^{\circ}] \times \C[\widehat{\Gr}_{4,N_R}^{\circ}]$ is a cluster algebra where each seed
is the disjoint union of a seed of each factor.} $\C[\widehat{\Gr}_{4,N_L}^{\circ}] \times \C[\widehat{\Gr}_{4,N_R}^{\circ}]$ to a sub-cluster algebra of $\C[\widehat{\Gr}_{4,n}^{\circ}]$. See \cite[Definition 3.1, Proposition 3.2]{Fraser} for the precise definition of a quasi-homomorphism.

\begin{theorem}\label{thm:promotion2}
Product promotion
$\Psi$ is a quasi-homomorphism of cluster
algebras. In particular, $\Psi$ maps a cluster variable (respectively, cluster)
of  $\C[\widehat{\Gr}_{4,N_L}^{\circ}] \times \C[\widehat{\Gr}_{4,N_R}^{\circ}]$, to a cluster variable (respectively, sub-cluster) of
$\C[\widehat{\Gr}_{4,n}^{\circ}]$, up to multiplication by Laurent monomials in $$\mathcal{T'}:=\{
\lr{a\,b\,c\,n},
\lr{a\,b\,c\,d},
\lr{b\,c\,d\,n},
\lr{a\,c\,d\,n}\}$$
\end{theorem}


Product promotion sends cluster variables to cluster variables up to scaling by a Laurent monomial factor. It will also be convenient to strip off this factor.

\begin{definition}\label{def:rPsi} Let $x$ be a cluster variable of $\C[\widehat{\Gr}_{4,N_L}^{\circ}]$ or $\C[\widehat{\Gr}_{4,N_R}^{\circ}]$.
 We define the \emph{rescaled product promotion} $\rPsi(x)$ of $x$ to be the cluster variable of $\Gr_{4,n}$ obtained from $\Psi(x)$ by removing the Laurent monomial\footnote{If $x= \lr{bcdn}$, then $\rPsi(x)=\Psi(x)= x$.} in $\mathcal{T'}$ (c.f. \Cref{thm:promotion2}).
\end{definition}

\begin{example}\label{ex:psi}
For $N_L$ and $N_R$ as in \cref{ex:bcw_tile}, we have
\[\Psi(\lr{1 \, 2 \, 4 \, 7})= \frac{\lr{1 \, 2 \, 4 \, 7}\lr{3\, 5\, 6\, 7}-\lr{1 \, 2 \, 3 \, 7}\lr{4\, 5\, 6\, 7} }{\lr{3 \, 4\, 6\, 7}}.\]
The remaining Pl\"uckers $\lr{1 \, 2 \, 3 \, 4}$, $\lr{1 \, 2 \, 3 \, 7}$, $\lr{1 \, 3 \, 4 \, 7}$, $\lr{2 \, 3 \, 4 \, 7}$ are fixed by $\Psi$.
We have $\rPsi(\lr{1 \, 2 \, 4 \, 7})=\lr{1 \, 2 \, 4 \, 7}\lr{3\, 5\, 6\, 7}-\lr{1 \, 2 \, 3 \, 7}\lr{4\, 5\, 6\, 7}$.  Note that this is a quadratic cluster variable in $\Gr_{4,7}$, obtained by mutating $\lr{2367}$ in $\Sigma_{4,7}$ of \cref{fig:gr_quiver}.
\end{example}

The fact that product promotion is a cluster quasi-homomorphism may be of independent interest
in the study of the cluster structure on $\Gr_{4,n}$.
Much of the work thus far on the cluster structure of the Grassmannian
has focused on cluster variables which are polynomials in Pl\"ucker coordinates
with low degree; by constrast, the cluster variables obtained by repeated application of $\rPsi$ can have arbitrarily high degree
in Pl\"ucker coordinates (e.g. see the \emph{chain polynomials} in \cite[Theorem 8.3]{Even-Zohar:2023del}).


Using rescaled product promotion and the operations in \cref{prop:pre_cluster},
we associate to each BCFW tile $\gt{S}$ a collection of compatible cluster variables $\Irr(S)$ for $\Gr_{4,n}$. We ultimately use these cluster variables to invert $\tilde{Z}$ on the tile, and to give a semi-algebraic description of each tile.

\begin{definition}[Cluster variables for BCFW tiles]\label{def:coord_cluster_vars}
Let $S \subset \Grk$ be a BCFW cell. We define the set of \emph{coordinate cluster variables} $\Irr(S)$ for $S$ recursively as follows:
\begin{itemize}
\item If $S =S_{L} \bowtie S_{R}$ with indices $B=(a, b, c, d, n)$, then
\begin{equation*}
    \Irr(S)=\rPsi_{B}(\Irr(S_L) \cup \Irr(S_R)) \cup \left\{ \lr{I}, I \in \binom{B}{4} \right\}.
\end{equation*}
\item If $S=\begin{cases}
\pre_i(S')& \\
\cyc(S')& \\
\refl(S')
\end{cases}$ then $\Irr(S)=\begin{cases}
\Irr(S')& \\
(\cyc^{-1})^*(\Irr(S'))& \\
\refl^*(\Irr(S'))
\end{cases}.$
\end{itemize}
For the base case $k=0$, we set $\Irr(S)=\emptyset$.
\end{definition}

For a BCFW cell $S$, $\Irr(S)$ depends on the sequence of operations in \cref{def:BCFW_cell} used to build $S$, but we will drop this dependence for brevity. Note that $\Irr(S)$ is a collection of compatible cluster variables for $\Gr_{4,n}$ \cite[Lemma 7.6]{Even-Zohar:2023del}.

\begin{example}\label{ex:coord_clust_vars} From \cref{ex:bcw_tile}, $S_{ex}=S_L \bowtie S_R$ with $B=(3,4,5,6,7)$. We have
 $\Irr(S_L)=\{\lr{I}, I \in \binom{B_L}{4} \}$ where $B_L=(1,2,3,4,7)$, and $\Irr(S_R)=\emptyset$. Then by \cref{ex:psi} the coordinate cluster variables $\Irr(S_{ex})$ are:
 $\lr{1234}$, $\lr{1237}$, $\lr{12  4  7}\lr{3 5 6 7}-\lr{1  2  3  7}\lr{4 5 6 7},\lr{1 3 4 7}$, $\lr{2 3 4 7}$ from
$\rPsi(\Irr(S_L))$ together with $\lr{3456}, \lr{3457}, \lr{3467}, \lr{3567}, \lr{4567}$.

\end{example}

Recall from \cref{def:tw_coords,def:functionary} the notion of twistor coordinates $\llrr{I}$ and functionaries. Given a cluster variable $x$ in $\Gr_{4,n}$, we let $x(Y)$ denote the functionary on $\Gr_{k,k+4}$ obtained by identifying Pl\"ucker coordinates $\lr{I}$ in $\Gr_{4,n}$ with twistor coordinates $\llrr{I}$ in $\Gr_{k,k+4}$.
Interpreting each cluster variable as a functionary in this way,
we describe each BCFW tile as the semialgebraic subset of
$\Gr_{k,k+4}$ where the coordinate cluster variables take on particular signs.

\begin{theorem}[Sign description of BCFW tiles]\label{th:sign_description}
	Let $\gt{S}$ be a BCFW tile. For each element $x$ of $\Irr(S)$, the functionary $x(Y)$ has a definite sign $s_x$ on the open tile $\gto{S}$ and
	\[\gto{S}= \{Y \in \Gr_{k,k+4}: s_x \, x(Y) >0 \text{ for all } x \in \Irr(S) \}.\]
\end{theorem}

\begin{example}
The open tile $\gto{S_{ex}}=\tilde{Z}(S_{ex})$, with $S_{ex}$ from \cref{ex:bcw_tile}, is the semialgebraic set in $\Gr_{2,6}$ described by: $x(Y)$ is negative for $x \in \{\lr{3567},\lr{3457},\lr{2347},\lr{3567}\}$ and $x(Y)$ is positive for all other $x \in \Irr(S_{ex})$ (computed in \cref{ex:coord_clust_vars}).
\end{example}
\cref{th:sign_description} is particularly relevant to \cref{thm:bcfw_tiles} as the ability to describe tiles using cluster variables follows from the fact that we can invert the amplituhedron map using cluster variables \cite[Theorem 7.7]{Even-Zohar:2023del}.

We now turn to facets of BCFW tiles (cf. \cref{def:facet2}); we will give functionaries which vanish on facets. The next result shows that each facet lies in the vanishing locus of a coordinate cluster variable (interpreting coordinate cluster variables as functionaries as above). 

\begin{theorem}[Cluster adjacency for BCFW tiles]
	\label{thm:clusteradjacency}
	Let $\gt{S}$ be a BCFW tile of
	$\Ank$. Each facet $\gt{S'}$ of $\gt{S}$ lies on a hypersurface cut out by a functionary $F_{S'}(\llrr{I})$ such that $F_{S'}(\lr{I})$ is in $\Irr(S)$. Thus $\{F_{S'}(\lr{I}): \gt{S'} \text{ a facet of }\gt{S}\}$ is a collection of compatible cluster variables of $\Gr_{4,n}$.
\end{theorem}

Finally, we explain how to construct BCFW tilings of $\Ank$ (cf. \cref{def:tiling}). This result shows that all ways of running the BCFW recurrence correspond to tilings of $\Ank$. \Cref{thm:bcfw_tiles,th:sign_description} are important ingredients to prove \cref{thm:BCFWtiling}.

In what follows, we use \cref{not:LR_cluster}, fix $n\geq k+4$, and define $b_{min}:=2$ if $k_L=0$ and otherwise $b_{min}:=k_L+3$.

\begin{definition}[BCFW collections]\label{def:BCFWoutput}

We say a collection $\mathcal{T}$ of $4k$-dimensional BCFW cells in
$\Gr_{k,n}^{\geq 0}$
	is a \emph{BCFW collection of cells} for
$\mathcal{A}_{n,k,4}$ if it has the following recursive form:
\begin{itemize}[noitemsep]
\item If $k=0$ or $k=n-4$, $\mathcal{T}$ is the single BCFW cell
		$\Gr_{k,n}^{>0}$.
\item If $\mathcal{T}=\{S\}$ is a BCFW collection of cells,
	so is $\{\refl S\}_{S \in \mathcal{T}}$ and $\{\cyc S\}_{S \in \mathcal{T}}$. \vspace{0.2em}
\item Otherwise $\mathcal{T} =
	\mathcal{T}_{pre} \bigsqcup_{b,k_L,k_R} \mathcal{T}_{k_L,k_R,b}$,
	where $k_L,k_R$ are as in \cref{not:LR_cluster}, $b$ ranges from
	$b_{min}$ to $n-3-k_R$, and
		\begin{itemize}
  \item
$\mathcal{T}_{pre}=\{\pre_{d}(S)\}_{S \in \mathcal{C}},$ where
			$\mathcal{C}$ is a BCFW collection of cells for $\mathcal{A}_{[n]\setminus\{d\},k, 4}$;
\item  $\mathcal{T}_{k_L,k_R,b}=\{S_L \bcfw S_R\}_{(S_L,S_R) \in \mathcal{C}_L \times \mathcal{C}_R}$
				where
				$\mathcal{C}_L$ and
				$\mathcal{C}_R$ are BCFW collections of cells for	$\mathcal{A}_{N_L,k_L,4}$
				and $\mathcal{A}_{N_R,k_R,4}$.
		\end{itemize}
		\end{itemize}
	\end{definition}

\begin{figure}
	\includegraphics[width=1.0\linewidth]{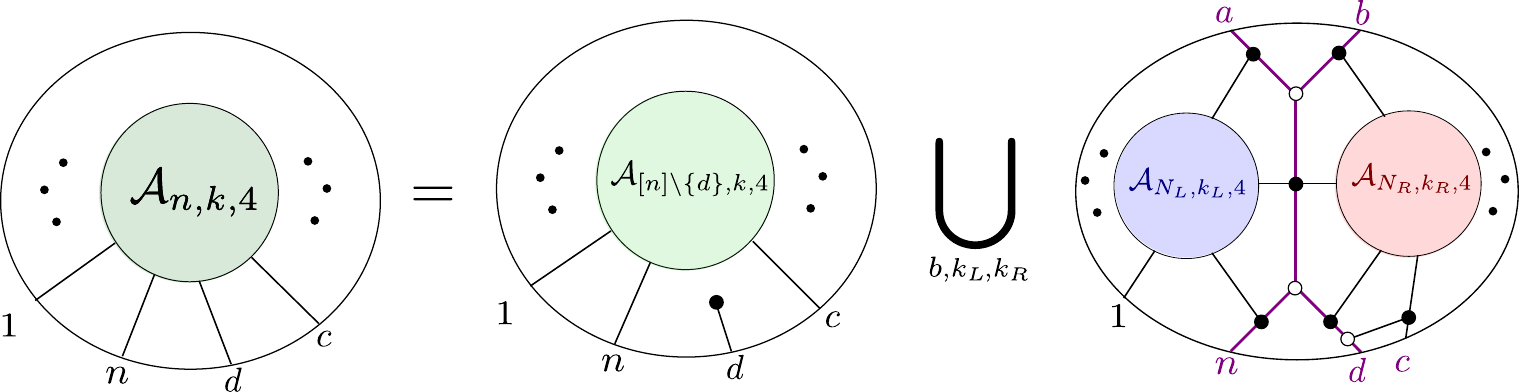}
	\caption{\label{fig:bcfw_tiling}
	BCFW tiling for $\mathcal{A}_{n,k,4}$. On the right: the first term is obtained by tiling $\mathcal{A}_{[n] \setminus \{d\},k,4}$ (from $\mathcal{T}_{pre}$); the second term is the union over $b,k_L,k_R$ as in \cref{def:BCFWoutput} of the collections of tiles obtained by tiling $\mathcal{A}_{N_L,k_L,4}$ and $\mathcal{A}_{N_R,k_R,4}$ (from $\mathcal{T}_{k_L,k_R,b}$).}
\end{figure}

The definition of BCFW collections comes directly from the BCFW recursion. Written graphically (see \cref{fig:bcfw_tiling}), the BCFW recursion expresses the scattering amplitude as a sum of terms indexed by cells in a BCFW collection.

\begin{theorem}[BCFW tilings]\label{thm:BCFWtiling}
Every BCFW collection of cells $\mathcal{T}=\{S \}$ as in \cref{def:BCFWoutput} gives rise to a tiling
	$\{Z_{S} \}_{S \in \mathcal{T}}$
	of the amplituhedron
$\Ank$, which we refer to as a \emph{BCFW tiling}.
\end{theorem}

  \cref{thm:BCFWtiling} generalizes the main result
of \cite{even2021amplituhedron}, which proved the same result for the \emph{standard} BCFW cells. \cref{thm:BCFWtiling} also proves
the main conjecture of \cite{arkani-hamed_trnka}.

Tiles which do not come from the BCFW recurrence are also expected to satisfy cluster adjacency, have a sign description in terms of cluster variables, and appear in tilings of $\Ank$. By leveraging the methods presented in this paper, we provided an example in \cite[Section 5]{Even-Zohar:2024nvw}.

\subsection*{Standard BCFW tiles}
 A BCFW cell is called \emph{standard} if it is obtained from the trivial cell $\Gr_{0, n}^{>0}$ by a series of BCFW products and applying $\pre_i$ at the penultimate index. The \emph{standard BCFW tiles} are the images of standard BCFW cells under the amplituhedron map. The collection of standard BCFW tiles for a fixed $n, k$ forms one tiling of the amplituhedron \cite{even2021amplituhedron}. We can sharpen a number of our results in the standard setting. Standard BCFW tiles are in bijection with \emph{chord diagrams}.

 \begin{definition}[Chord diagram \cite{even2021amplituhedron}]\label{def:cd}
Let $k,n \in \mathbb{N}$. A~\emph{chord diagram} $D \in\mathcal{CD}_{n,k}$ is a set of $k$~quadruples, called \emph{chords}, in the \emph{marker set} $\{1,\dots,n\}$
	$$ D \;=\; \{(a_1,b_1,c_1,d_1),\dots,(a_k,b_k,c_k,d_k)\}$$ such that every
	chord $D_i=(a_i,b_i,c_i,d_i) \in D$ satisfies
$ 1 \;\leq\; a_i \;<\; b_i=a_i+1 \;<\; c_i \;<\; d_i=c_i+1 \;\leq\; n-1 $
and \emph{no} two chords $D_i,D_j \in D$ satisfy
$ a_i \;=\; a_j$ or $a_i \;<\; a_j \;<\; c_i \;<\; c_j.$
\end{definition}
\begin{figure}
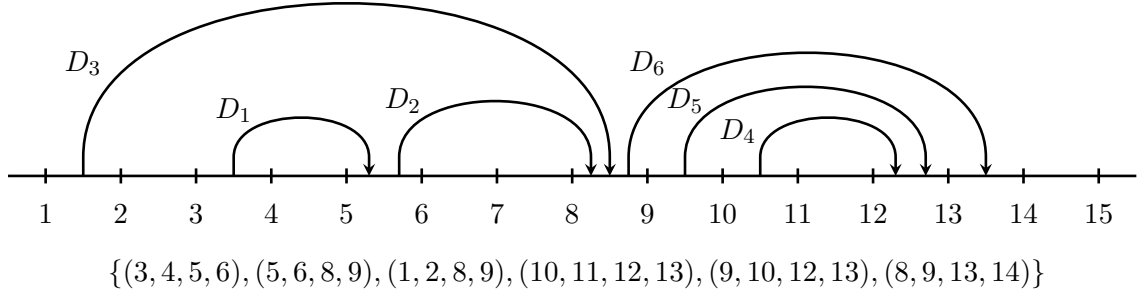

\begin{center}
\tikz[line width=1,scale=1]{
\draw (0.5,0) -- (15.5,0);
\foreach \i in {1,2,...,15}{
\def\x{\i}
\draw (\x,-0.1)--(\x,+0.1);
\node at (\x,-0.5) {\i};}
\foreach \i/\j in {1/8, 3/4.8, 5.2/7.75, 8.25/13, 9/12.2, 10/11.8}{
\def\x{\i+0.5}
\def\y{\j+0.5}
\draw[line width=1,-stealth] (\x,0) -- (\x,0.25) to[in=90,out=90] (\y,0.25) -- (\y,0);
}
\node at(1.5,1.5) {{$D_3$}};
\node at(3.5,0.875) {{$D_1$}};
\node at(5.75,1) {{$D_2$}};
\node at(9,1.5) {{$D_6$}};
\node at(9.5,1) {{$D_5$}};
\node at(10.2,0.6) {{$D_4$}};
}
\vspace{0.5em}

{$\{(3,4,5,6),(5,6,8,9),(1,2,8,9),(10,11,12,13),(9,10,12,13),(8,9,13,14)\}$}
\end{center}
\caption{
A chord diagram $D$ with $k=6$ chords and $n=15$ markers.
}
\label{cd-example}
\end{figure}

See \cref{cd-example} for an example. Notice that the chords are noncrossing and no two chords start in the same place. For the details of the correspondence between standard BCFW tiles and chord diagrams, see \cite[Definition 6.12]{Even-Zohar:2023del}; each chord $D_i$ corresponds to a BCFW product with $B_i=(a_i,b_i,c_i,d_i,n)$. Recall that the corresponding product promotion adds $\{\lr{I}: I \in \binom{B_i}{4}\}$ to the coordinate cluster variables. In what follows, we denote by $\ralpha_i, \rbeta_i, \rgamma_i, \rdelta_i, \repsilon_i$ the coordinate cluster variables obtained from $\lr{b_i c_i d_i n},\lr{a_i c_i d_i n},\lr{a_i b_i d_i n},\lr{a_i b_i c_i n},\lr{a_i b_i c_i d_i}$, respectively, by applying some product promotions.

A chord $D_i$ is a \emph{child} of the chord immediately above it, which is its \emph{parent}. For example, $D_4$ has parent $D_5$ in \cref{cd-example}. Two chords are \emph{same-end}
if their ends occur in a common segment $(e,e+1)$ (e.g. $D_2, D_3$); are \emph{head-to-tail} if the first ends in the segment where the second starts (e.g. $D_1, D_3$); and are \emph{sticky} if their
starts lie in consecutive segments $(s,s+1)$ and~$(s+1,s+2)$ (e.g. $D_5, D_6$).

Given a chord diagram $D$ corresponding to a standard BCFW cell $S_D$, we give explicit formulas for the coordinate cluster variables $\xx(S_D)$ and their signs on the BCFW tile $Z_D$ in \cite[Theorem 8.4, Proposition 8.10]{Even-Zohar:2023del}. We also determine which coordinate cluster variables cut out facets of $Z_D$.

\begin{theorem}
    Let $\rzeta_i$ be a coordinate cluster variable for $Z_D$, considered as a functionary. A facet of $Z_D$ lies on the hypersurface $\{\rzeta_i=0\}$ if and only if:
    \begin{itemize}[leftmargin=20pt]
        \item $\rzeta_i = \ralpha_i$ and $D_i$ does not have a sticky child, or
        \item $\rzeta_i = \rbeta_i$ and $D_i$ does not start where another chord ends or have a sticky same-end parent, or
        \item $\rzeta_i = \rgamma_i$, or
        \item $\rzeta_i \in \{\rdelta_i, \repsilon_i\}$ and $D_i$ does not have a same-end child.
\end{itemize}
\end{theorem}

We call $\rzeta_i \in \xx(S_D)$ \emph{frozen} if a facet of $Z_D$ lies in its zero locus; otherwise we call $\czeta_i$ \emph{mutable}. Recall that the coordinate cluster variables $\xx(S_D)$ are compatible cluster variables for $\Gr_{4,n}$, and so appear together in some seed $\Sigma(D)$. We additionally determine all arrows in the quiver for $\Sigma(D)$ which involve the mutable variables of $\xx(S_D)$.


\begin{theorem}\label{thm:quiver}
Let $D \in \mathcal{C}\mathcal{D}_{n,k}$, let $S_D$ be the standard BCFW cell, and let $\Sigma_D=(\xx, Q)$ be a seed for $\Gr_{4,n}$ whose cluster contains $\xx(S_D)$. Then the arrows of $Q$ which involve the mutable variables of $\xx(S_D)$ are as follows:
For each chord $D_i$,
check if it forms a configuration in the table below, and if so, draw the indicated arrows\footnote{In the third configuration, if $D_i$ and $D_j$ end in the same place, then the dotted arrow from $\ralpha_i$ to $\repsilon_i$ appears, along with the arrows from the second and third configuration.}.
\vspace{0.5em}
\begin{center}
\begin{tabular}{|c|c|c|}
\hline
\tikz[line width=1,scale=1]{
\def\r{1}
\draw (1*\r,0) -- (1.5*\r,0);
\draw[dashed] (1.5*\r,0) -- (2.75*\r,0);
\draw (2.75*\r,0) -- (3.5*\r,0);
\draw[dashed] (3.5*\r,0) -- (5*\r,0);
\draw (5*\r,0) -- (5.5*\r,0);
\foreach \i in {6,7}{
\def\x{\i/2*\r}
\draw (\x,-0.1)--(\x,+0.1);
}				
\foreach \i/\j in {2/5.833,6.166/10}{
\def\x{\i/2*\r+0.25*\r}
\def\y{\j/2*\r+0.25*\r}
\draw[line width=1,-stealth] (\x,0) -- (\x,0.15) to[in=90,out=90] (\y,0.15) -- (\y,0);}
\node at(2.2*\r,1*\r) {$j$};
\node at(4.3*\r,1*\r) {$i$};
}&
\tikz[line width=1,scale=1]{
\draw (0.75,0) -- (1.25,0);
\draw[dashed] (1.25,0) -- (2.25,0);
\draw (2.25,0) -- (2.75,0);
\draw[dashed] (2.75,0) -- (4,0);
\draw (4,0) -- (5,0);
\foreach \i in {8.5,9.5}{
\def\x{\i/2}
\draw (\x,-0.1)--(\x,+0.1);
}
\foreach \i/\j in {1/4.5+0.066, 2.5/4.5-0.066}{
\def\x{\i}
\def\y{\j}
\draw[line width=1,-stealth] (\x,0) -- (\x,0.15) to[in=90,out=90] (\y,0.15) -- (\y,0);}
\node at(2.4,0.6) {$j$};
\node at(2.7,1.5) {$i$};
}&
\tikz[line width=1,scale=1]{
\draw (0.5,0) -- (1.75,0);
\draw[dashed] (1.75,0) -- (3.25,0);
\draw (3.25,0) -- (3.75,0);
\draw[dashed] (3.75,0) -- (4.75,0);
\draw (4.75,0) -- (5.25,0);
\foreach \i in {1.5,2.5,3.5}{
\def\x{\i/2}
\draw (\x,-0.1)--(\x,+0.1);
}
\foreach \i/\j in {1/5, 1.5/3.5}{
\def\x{\i}
\def\y{\j}
\draw[line width=1,-stealth] (\x,0) -- (\x,0.15) to[in=90,out=90] (\y,0.15) -- (\y,0);}
\node at(3,1) {$j$};
\node at(2,1.5) {$i$};
}\\
head-to-tail left sibling $D_j$
& same-end child $D_j$
& sticky child $D_j$
\\[0.5em]
\hline
\tikz[line width=0.75,scale=1]{
\node        (c1) at (4,5.5) {$\rgamma_j$};
\node        (d1) at (6,5.5) {$\rdelta_j$};
\node        (a2) at (4,4) {$\ralpha_i$};
\node (b2) at (6,4) {$\rbeta_i$};
\path[very thick,->] (b2) edge (d1);
\path[very thick,->] (c1) edge (b2);
\path[very thick,->] (b2) edge (a2);
}
&
\tikz[line width=0.8,scale=1]{
\node        (c5) at (12,2.5) {$\rgamma_i$};
\node       (d5) at (13.5,2.5) {$\rdelta_i$};
\node        (e5) at (15,2.5) {$\repsilon_i$};
\node        (c4) at (12,1) {$\rgamma_j$};
\node        (d4) at (13.5,1) {$\rdelta_j$};
\node        (e4) at (15,1) {$\repsilon_j$};
\path[very thick,->] (e5) edge (d5);
\path[very thick,->] (d5) edge (c5);
\path[very thick,->] (c4) edge (d5);
\path[very thick,->] (d4) edge (e5);
\path[very thick,->] (d5) edge (d4);
\path[very thick,->] (e5) edge (e4);
}
&
\tikz[line width=0.75,scale=1]{
\node    (a5) at (8.5,2) {$\ralpha_i$};
\node (b5) at (10,2) {$\rbeta_i$};
\node (e5) at (6.5,2) {$\repsilon_i$};
\node (a4) at (10,0.5) {$\ralpha_j$};
\node (e4) at (6.5,0.5) {$\repsilon_j$};
\node[gray] (text) at (7.5,2.5) {if same-end};
\path[very thick,->,dotted] (a5) edge (e5);
\path[very thick,->] (e4) edge (a5);
\path[very thick,->] (b5) edge (a5);
\path[very thick,->] (a5) edge (a4);
}
\\[0.5em]
\hline
\end{tabular}
\end{center}
\vspace{0.5em}

\end{theorem}

\begin{figure}
    \centering
    \includegraphics[width=0.5\linewidth]{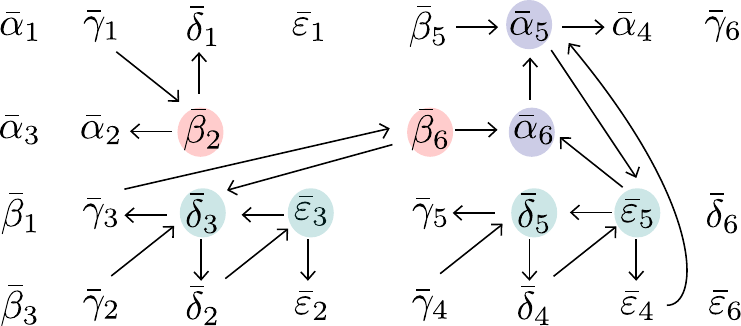}
    \caption{The quiver $Q_D$ for the chord diagram in \cref{cd-example}. The highlighted variables are the mutable variables, which do not cut out facets of the tile $Z_D$.}
    \label{fig:quiver-example}
\end{figure}

\noindent
\textbf{Acknowledgements.}
TL is supported by SNSF grant Dynamical Systems, grant no.~188535.
MP is supported by Harvard CMSA and at IAS by the grant number DE-SC0009988.
MSB is supported by the National Science Foundation under Award No.~DMS-2103282.
RT (incumbent of the Lillian and George Lyttle Career Development Chair) was supported by the ISF grant No.~335/19 and 1729/23.
LW is supported by the National Science Foundation under Award No.
DMS-1854316 and DMS-2152991. Any opinions, findings, and conclusions or recommendations expressed in this material are
those of the author(s) and do not necessarily reflect the views of the National Science
Foundation.

\bibliographystyle{alpha}
\bibliography{one-column}

\end{document}